\documentclass[12pt]{article}
\usepackage{amsmath}
\usepackage{amsthm}
\usepackage{amssymb}

\newtheorem{thm}{Theorem}

\newtheorem{prop}{Proposition}
\newtheorem{lemma}{Lemma}
\newtheorem{definition}{Definition}

\begin{document}

\vspace*{30px}

\begin{center}\Large
\textbf{Hadamard Structures with Associated Automorphisms}
\bigskip\large

Ivica Martinjak\\
University of Zagreb, Faculty of Science\\
Croatia


\end{center}





\begin{abstract}
In this paper we present new Hadamard matrices and related combinatorial structures. In particular, it is constructed 5202 inequivalent Hadamard matrices of order 36  as well as 180538 Hadamard symmetric designs with 35 points in addition to those structures that admit an automorphism of order 3. Consequently, there are at least 272116 Hadamard 3-designs with 36 points and 70 lines. We found that all Hadamard matrices constructed here are equivalent to a regular Hadamard matrix. This fact contributes to the conjecture that Hadamard matrices of order 36, and possibly those of order $4m^2$, are regular.

\end{abstract}






\section{Introduction}
\label{Intro}

A {\it Hadamard matrix} $H$ of order $n$ is a square matrix with entries $\pm1$ satisfying  $HH^t=nI_n$, where $I_n$ is the identity matrix of order $n$. Having determinant equal to $\pm n^{n/2}$, these matrices are solution to the Hadamard's Maximum Determinant Problem. It can be easily shown that if $n$ is the order of $H$ then $n=1,2$ or $n \equiv 0 \pmod 4$. However, the conjecture that for every natural number $n$ divisible by $4$ there exists a Hadamard matrix of order $n$ is still an open problem. 

A Hadamard matrix with all-ones first row and column is called $normalized$. A Hadamard matrix $H$ of order $n$ is skew if $H^t +H =2I_n$. One of a few notable properties of Hadamard matrices is also {\it regularity}. A matrix $H$ having constant row and column sum is called regular. Moreover, a Hadamard matrix of order $4m^2$ is regular if every row and column of $H$ contains a constant number $(2m^2-m$ or $2m^2+m)$ of $+1$'s. The necessary condition for regularity is that the order of $H$ be a perfect square. Having in mind that Hadamard matrices can be interpreted as 2-designs, this means that Menon-type series of designs possibly lead to the regular Hadamard matrices.

The following operations on Hadamard matrices preserve the Hadamard property: {\it i)} permuting rows or columns, {\it ii)} negating rows or columns, {\it iii)} transposition.

Two Hadamard matrices $H_a$ and $H_b$ are called {\it equivalent} if one can be obtained from the other by operation of types {\it i)} and {\it ii)}. In other words, matrices $H_a$ and $H_b$ are equivalent if 
\begin{eqnarray}
H_b=P^{-1}H_aQ, 
\end{eqnarray}
where $P$ and $Q$ are monomial matrices. Up to equivalence, there is a unique Hadamard matrix of orders up to 12. There are 5, 3, 60 and 487  inequivalent matrices of the following orders up to 28. The highest order for which complete classification is known is 32, which is recently done by Kharaghani and Tayfeh-Rezaie \cite{Khar}.

Thus, the next order for which a complete classification is not known is 36. In addition, Hadamard matrices of this order are of particular interest because of they rich relation with the other combinatorial and algebraic structures. There are 11 matrices of this order admitting an automorphism of order 17, which are found by V. Tonchev \cite{Tonc}. While it is known that symmetric Bush-type matrices of this order does not exist, Z. Janko constructed a nonsymmetric Bush-type matrix  \cite{Jank}. The most comprehensive study of matrices of this order is done in \cite{BFW}. In that work $38332$ matrices is constructed, that arise from an action of automorphisms of order 2, 3 5, 7 and 17. Remarkably, all of these matrices are regular.

\section{Hadamard 3-designs}
\label{3-des}

Recall that an {\it incidence structure} is a triple ${\cal I}=({\cal P}, {\cal L}, I)$ where ${\cal P}=\{p_1,...,p_v\}$ is a set of {\it points}, the elements of a set ${\cal L}=\{L_1,...,L_b\}$ are called {\it lines} or {\it blocks} and $I \subseteq {\cal P} \times {\cal L}$ is an {\it incidence relation}. The incidence structure having every $t$ points on the same number $\lambda$ of lines and having the property that every line is incident with the same number $k$ of points is of particular interest. Such structure is called $t$-design and it is uniquely determined by the 4-tuple of parameters $t$-$(v,k,\lambda)$. Every point $p \in {\cal P}$ of a $t$-design is incident with the same number of lines $r$: 
\begin{eqnarray}  \label{r_defi}
r \binom{k-1}{t-1}=\lambda \binom{v-1}{t-1}.
\end{eqnarray} 
{\it Symmetric designs} are special class of $t$-designs for $t=2$, having the number of points $v$ equal to the number of lines $b$. More details on such structures one can find in books and handbooks on combinatorial design theory \cite{HuPu}, \cite{ITT}.

Let ${\cal D}=({\cal P},{\cal B},I)$ and ${\cal D'}=({\cal P'},{\cal B'},I')$ be two $t$-designs. A bijection ${\phi}:{\cal P} \cup {\cal B} \to {\cal P'} \cup {\cal B'}$ is an isomorphism if $\phi$ maps points onto points and blocks onto blocks, and $(p,B)\in I \Leftrightarrow ({\phi}(p),{\phi}(B))\in I', \enspace \forall p \in {\cal P}, \enspace \forall B \in {\cal B}.$ Two $t$-designs $\cal D$ and $\cal D'$ are {\it isomorphic}, $\cal D \approx \cal D'$, if there exists an isomorphism of $\cal D$ on $\cal D'$. An isomorphism of $t$-design $\cal D$ to itself we call an {\it automorphism} of $\cal D$. It is known that a set of all automorphisms of $\cal D$ form a group, which is called a {\it a full automorphisms group} and it is denoted by $Aut({\cal D})$. An {\it automorphism group order} $|Aut({\cal D})|$ is the cardinality of a full automorphism group of $\cal D$.

It can be said that $3$-designs are even more regular structure than $2$-designs since within $3$-designs every triple of {\it incidence matrix} rows intersect in $\lambda$ points. Clearly, every $t$-design, $t\ge 2$ is in the same time a $2$-design. More precisely, for every $t$-$(v,k,\lambda)$ design there exists a $s$-$(v,k,\lambda_s)$ design, $0 \leq s \leq t$, where
\begin{eqnarray}
\lambda_s = \lambda \frac{(v-s)(v-s-1)...(v-t+1)}{(k-s)(k-s-1)...(k-t+1)}.
\end{eqnarray}

A Hadamard 3-design ${\cal D}$ with parameters $3$-$(4a,2a,a-1)$ exists if and only if there exist a Hadamard matrix of order $n=4a$. Namely, starting from a matrix $H$ of order $n=4a$, one can construct an incidence matrix $A$ of ${\cal D}$. When fixing a row of $H$, then each of the other $n-1$ rows $s$ defines two rows of the matrix $A$: {\it i)} the first row is obtained by putting '1' in any position in which chosen row mach with the row $s$; {\it ii)} the second row is the complement of the first one. Hadamard $3$-designs constructed from the same Hadamard matrix can be non-isomorphic.

There is a unique Hadamard 3-design for $a=1,2,3$. It is also known that there are 5, 3 and 130 structures when $a$ is equal to 4, 5 and 6, respectively. It is worth mentioning that some Hadamard 3-designs of Paley type are constructed in \cite{BBI}.

Let $H$ be a normalized Hadamard matrix of order $n=4a$. Then removing the first row and column of $H$ and replacing the -1's with 0 in $H$ gives the incidence matrix of a symmetric $2$-$(4a-1,2a-1,a-1)$ design. The opposite procedure holds as well. Designs with these parameters are called Hadamard 2-designs.

In this paper we aim at constructing Hadamard matrices, Hadamard 2-designs and Hadamard 3-designs for $a=m^2=9$, that admit a {\it tactical decomposition}. In addition, for the purpose to reach new structures we also perform constructions ignoring a group action in the second step of the procedure. Having known that it is not known any exception of the fact that already found Hadamard matrices of order 36 are regular, our objective is also to see whether new matrices follow this property.

\section{Method of construction}

\begin{definition}
Let ${\cal D}=({\cal P}, {\cal B}, I)$ be a $t$-$design$, and the decomposition 
$${\cal P} = {\cal P}_1 \sqcup {\cal P}_2 \sqcup
\cdots \sqcup {\cal P}_m \mbox{ and } {\cal B} = {\cal B}_1 \sqcup {\cal
B}_2 \sqcup \cdots \sqcup {\cal B}_n\, $$ 
be a partition of ${\cal P}$ and ${\cal B}$.
Furthermore, let every point in a set ${\cal P}_i, \enspace i=1,\dots,m$ be incident with the same number of blocks in a set ${\cal B}_j, \enspace j=1,\dots,n$; and each blocks from ${\cal B}_i, \enspace i=1,\dots,n$ incident with the same number of points in a set ${\cal P}_j, \enspace j=1,\dots,m$. Then the decomposition is {\it tactical}. 
\end{definition}

If conditions of this definition are fulfilled, an incidence matrix $M=[m_{ij}]$ of a design $\cal D$ has a significant characteristic, as follows. Obtained tactical decomposition then naturally divide matrix $M$ on submatrices $A_{ij}$, $1 \le and \le m$, $1 \le j \le n$ where particular submatrix $A_{ij}$ has dimension $\vert {\cal P}_i \vert \times \vert {\cal B}_j \vert$;
$$
M=
\left[\begin{array}{cccc}
A_{11} & A_{12} & ... & A_{1n} \\
A_{21} & A_{22} & ... & A_{2n} \\
\vdots & \vdots &     & \vdots \\
A_{m1} & A_{m2} & ... & A_{mn}\end{array}\right].
$$
These submatrices have the same number of 1's in every row as well in every column. Conversely, if submatrices possess this characteristics there holds conditions from the previous definition. If a decomposition of an incidence structure ${\cal I}=({\cal P}, {\cal B}, I)$ is tactical, than we can define the next coefficients
\begin{eqnarray*}
\rho_{ij} & = & |\{ B \in {\cal B}_j \mid (p,B)\in I, \enspace p \in {\cal P}_i \}|, \\ 
\kappa_{ij} & = & |\{ p \in {\cal P}_i \mid (p,B)\in I, \enspace B \in {\cal B}_j \} \vert.
\end{eqnarray*}   
These coefficients represent the number of blocks from a set ${\cal B}_j$ which are incident with a particular point from a set ${\cal P}_i$, and the number of point in a set ${\cal P}_i$ which contained in a particular block of a set ${\cal B}_j$, where $i=1,\dots,m$ and $j=1,\dots,n$. Matrices $[\rho_{ij}]$ and $[\kappa_{ij}]$ are called {\it matrices of tactical decomposition of a design} $\cal D$.

\begin{lemma}
Let ${\cal D}=({\cal P}, {\cal B}, I)$ be a symmetric $(v,k,\lambda)$ design, and let 
$${\cal P} = {\cal P}_1 \sqcup {\cal P}_2 \sqcup
\cdots \sqcup {\cal P}_m \mbox{,    } \enspace {\cal B} = {\cal B}_1 \sqcup {\cal
B}_2 \sqcup \cdots \sqcup {\cal B}_n\, $$ 
be its tactical decomposition. 
Then the following equations hold for coefficients of tactical decomposition matrices (TDM's). 
\begin{eqnarray}
\sum\limits_{j=1}^m \rho_{ij} & = & k, \enspace \forall i \label{relk}\\
\sum_{j=1}^n \frac{|{\mathcal P}_l|}{|{\mathcal B}_j|} \rho_{ij}
\rho_{lj} &=& \lambda \cdot |{\mathcal P}_l| +
\delta_{il}(k-\lambda) \label{skalprod2}
\end{eqnarray}
where $\delta_{il}$ is the Kronecker symbol.
\end{lemma}

Let ${\cal D}$ be an incidence structure and $G$ a subgroup of $Aut({\cal D})$. It is well known that then $G$ forms a tactical decomposition of ${\cal D}$ (for details and prove see \cite{KaOs}). This fact allow us a two-step construction procedure that finds firstly all non-isomorphic solutions of the equation system (\ref{relk}) and (\ref{skalprod2}), as it is applied in \cite{CEPU} and \cite{Krc}. The second step extends these matrices by determining exactly the incidences between points and lines \cite{JTT}. These means that combination of cyclic matrices should be checked in respect to defining properties of a design. Recall that a binary $v \times b$ matrix is an incidence matrix of a $t$-design  if and only if {\it i)} every row has $r$ 1's, {\it ii)} every column has $k$ 1's and {\it iii)} the scalar product of every $t$ rows is $\lambda$. These properties follows from the defining properties of a design and it also can be shown that the second and third property imply the first one. 

In order to reach new enumeration results and to get a better insight into the nature of Hadamard structures, in particular those for $m=3$, we also use a variation of this approach. Namely, in the second step of our construction method we ignore a group action. This means that anti-cyclic matrix were also included in possibly constructed design. However, in this case for every entry $\rho_{ij}=1$ or 2, the number of possibilities for substituting it with a 0-1 matrix was doubled; leading to the significantly longer exhaustive search \cite{MaPa}. Throughout the paper, we will use abbreviation Cyc for the first type of a construction and ACyc for the second one.

\section{Constructions and partial classification of Hadamard structures for $m=3$}
\label{Results}

\subsection{Results for Hadamard BIBD's for 2-(35,17,8)}
An automorphism of order $3$ acts on a $2$-$(35,17,8)$ design fixing $F \in \{2,5,8\}$ points and lines. This result is in line with facts that the design is symmetric and orbits are of length 1 or 3. For the most complex case when $F=2$, we obtained 626 tactical decomposition matrices. Furthermore, there are 16 TDM's for $F=5$ and two matrices when the automorphism fixes $8$ points and blocks. Once these matrices are obtained we perform the standard tactical decomposition method as it is described above, whose second phase results with incidence matrices. Finally, non-isomorphic structures among these matrices are isolated \cite{McK2}. We obtain 63635, 3698 and 14692 non-isomorphic designs for 2, 5 and 8 fixed points, respectively.

Having in mind that designs with the different number of fixed points and lines can be isomorphic, at least we test isomorphisms of these three sets of incidence matrices. This proves the next theorem.

\begin{thm}
There are exactly $81973$ symmetric designs with parameters $(35,17,8)$ admitting an action of an automorphism of order $3$.
\end{thm}

The most numerous among them are designs having $|Aut({\cal D})|=3$ while the most symmetric design has the automorphism group order of 40320. It follows the complete list of automorphism group orders appearing, together with the next frequencies: \\
3 (79704), 6 (1916), 9 (30), 12 (156), 18 (43), 21 (1), 24 (49), 36 (13), 48 (4), 60 (1), 72 (6), 96 (4), 144 (2), 168 (1), 192 (1), 288 (1), 384 (1), 420 (1), 720 (1), 1152 (1), 40320 (1).

In order to construct new structures we performed the ACyc construction method, where all permutation matrices of order 3 are taken into account. This time we get 89686,  17627 and 155550 non-isomorphic designs for using TDM's with 2, 5 and 8 fixed point, respectively. These enumeration results are summarized in the next proposition.

\begin{prop}
There are at least  $262511$ symmetric designs with parameters $(35,17,8)$, where $169002$ of them do not admit any non-trivial automorphisms.
\end{prop}

Thus, we were able to reach structures without any non-trivial automorphism and such structures dominate among constructed designs. More precisely, the frequencies of the already presented automorphism group orders remain the same in this construction. In addition to automorphism group orders appearing in Cyc construction, here we have the next values for $|Aut({\cal D})|$: 1 (169002), 2 (10776), 4 (663), 8 (102), 10 (2), 16 (18), 32 (8), 64 (3).

\subsection{Results for Hadamard matrices of order 36}
Once having $2$-$(35,17,8)$ designs constructed, we transform them into Hadamard matrices adding all-1 first row and column and changing 0's to -1's. Test on equivalence is performed using the fact that two Hadamard matrices are equivalent if certain related graphs are isomorphic \cite{McK1}. Detailed results are presented in Table \ref{Table1} where the number of inequivalent matrices corresponding to designs arising by Cyc construction is denoted by Cyc*; and we use an analogue notation for the another type of a construction. 

In total 12833 Hadamard matrices of order 36 is constructed, 7631 of them corresponding to the Hadamard 2-designs with 35 points admitting an automorphism group action. In other words there are 5202 matrices arising from designs that do not admit an action of an automorphism of order 3.

\begin{prop} There are at least $12833$ inequivalent Hadamard matrices of order $36$, where $7631$ out of them correspond to the Hadamard $2$-designs with parameters $2$-$(35,17,8)$ that admit an action of an automorphism of order $3$. All of these matrices are regular. 
\end{prop}

Being aware that all already known Hadamard matrices of this order are regular, we were intrigued to check weather the same is case for matrices constructed in this work. For this purpose we developed a custom algorithm that a Hadamard matrix transform into equivalent regular matrix if this is possible. Impressively, this operation has shown that all matrices constructed from 2-designs that do not admit an automorphism of order 3 are regular. In other words, for every of these matrices we were able to find an equivalent regular Hadamard matrix. Having in mind that all matrices constructed in \cite{BFW} are regular our result means that all currently known Hadamard matrices of order 36 are regular; which additionally contributes to the conjecture that all Hadamard matrices of order 36 are regular. As an illustration, the following regular Hadamard matrix $H_1$ of order 36 arise from the Hadamard 2-design with automorphism group order of 40320 (- means -1).

\begin{scriptsize}
\[ H_{1}= \left( \begin{tabular}{c@{\hspace{0.5em}}c@{\hspace{0.5em}}c@{\hspace{0.5em}}c@{\hspace{0.5em}}c@{\hspace{0.5em}}c@{\hspace{0.5em}}c@{\hspace{0.5em}}c@{\hspace{0.5em}}c@{\hspace{0.5em}}c@{\hspace{0.5em}}c@{\hspace{0.5em}}c@{\hspace{0.5em}}c@{\hspace{0.5em}}c@{\hspace{0.5em}}c@{\hspace{0.5em}}c@{\hspace{0.5em}}c@{\hspace{0.5em}}c@{\hspace{0.5em}}c@{\hspace{0.5em}}c@{\hspace{0.5em}}c@{\hspace{0.5em}}c@{\hspace{0.5em}}c@{\hspace{0.5em}}c@{\hspace{0.5em}}c@{\hspace{0.5em}}c@{\hspace{0.5em}}c@{\hspace{0.5em}}c@{\hspace{0.5em}}c@{\hspace{0.5em}}c@{\hspace{0.5em}}c@{\hspace{0.5em}}c@{\hspace{0.5em}}c@{\hspace{0.5em}}c@{\hspace{0.5em}}c@{\hspace{0.5em}}c@{\hspace{0.5em}}} 
 1& 1& 1& -& -& -& -& -& -& -& -& 1& -& 1& -& 1& -& -& -& 1& -& 1& -& -& -& -& 1& -& 1& 1& 1& -& 1& 1& 1& -\\
 1& 1& 1& -& -& -& -& -& -& -& -& 1& -& 1& -& 1& -& -& 1& -& 1& -& 1& 1& 1& 1& -& 1& -& -& -& 1& -& -& -& 1\\
 1& 1& 1& -& -& -& -& -& -& 1& 1& -& 1& -& 1& -& 1& 1& -& 1& -& 1& -& -& -& -& 1& 1& -& -& -& 1& -& -& -& 1\\
 -& -& -& 1& 1& 1& -& -& -& 1& 1& -& -& 1& -& 1& -& -& 1& -& 1& 1& -& -& -& -& 1& 1& -& -& -& 1& -& 1& 1& -\\
 -& -& -& 1& 1& 1& -& -& -& -& -& 1& 1& -& 1& 1& -& -& -& 1& -& -& 1& 1& -& -& 1& 1& -& -& 1& -& 1& -& -& 1\\
 -& -& -& 1& 1& 1& -& -& -& -& -& 1& -& 1& -& -& 1& 1& -& 1& -& 1& -& -& 1& 1& -& -& 1& 1& -& 1& -& -& -& 1\\
 -& -& -& -& -& -& 1& 1& 1& 1& -& 1& 1& 1& -& -& -& -& 1& 1& -& -& -& -& 1& -& 1& 1& -& 1& -& 1& 1& -& -& -\\
 -& -& -& -& -& -& 1& 1& 1& -& 1& 1& -& -& -& 1& 1& -& -& -& -& 1& 1& -& -& 1& 1& -& -& -& 1& 1& -& 1& -& 1\\
 -& -& -& -& -& -& 1& 1& 1& -& -& -& -& 1& 1& 1& -& 1& -& 1& 1& 1& -& 1& -& -& -& 1& 1& -& -& -& -& -& 1& 1\\
 -& -& 1& 1& -& -& 1& -& -& 1& 1& -& 1& 1& -& -& -& -& 1& 1& -& 1& 1& 1& -& 1& -& -& 1& -& 1& -& -& -& -& -\\
 -& -& 1& 1& -& -& -& 1& -& 1& 1& -& -& -& -& 1& 1& -& -& -& -& -& -& 1& 1& -& -& 1& 1& 1& -& -& 1& 1& -& 1\\
 1& 1& -& -& 1& 1& 1& 1& -& -& -& 1& 1& -& -& -& 1& -& 1& -& -& 1& -& 1& -& -& -& 1& 1& -& -& -& -& 1& -& -\\
 -& -& 1& -& 1& -& 1& -& -& 1& -& 1& 1& -& 1& -& -& -& -& -& 1& -& -& -& -& 1& -& -& 1& -& -& 1& 1& 1& 1& 1\\
 1& 1& -& 1& -& 1& 1& -& 1& 1& -& -& -& 1& -& -& -& 1& -& -& -& -& -& 1& -& 1& 1& -& -& -& -& -& 1& 1& -& 1\\
 -& -& 1& -& 1& -& -& -& 1& -& -& -& 1& -& 1& 1& -& 1& 1& -& -& 1& -& 1& 1& 1& 1& -& -& 1& -& -& -& 1& -& -\\
 1& 1& -& 1& 1& -& -& 1& 1& -& 1& -& -& -& 1& 1& -& -& 1& 1& -& -& -& -& -& 1& -& -& 1& -& -& 1& 1& -& -& -\\
 -& -& 1& -& -& 1& -& 1& -& -& 1& 1& -& -& -& -& 1& 1& 1& 1& 1& -& -& 1& -& 1& 1& -& -& -& -& -& 1& -& 1& -\\
 -& -& 1& -& -& 1& -& -& 1& -& -& -& -& 1& 1& -& 1& 1& 1& -& -& -& 1& -& -& -& -& 1& 1& -& 1& 1& 1& 1& -& -\\
 -& 1& -& 1& -& -& 1& -& -& 1& -& 1& -& -& 1& 1& 1& 1& 1& -& 1& -& -& -& 1& -& 1& -& 1& -& 1& -& -& -& -& -\\
 1& -& 1& -& 1& 1& 1& -& 1& 1& -& -& -& -& -& 1& 1& -& -& 1& -& -& -& 1& 1& -& -& -& -& -& 1& 1& -& -& 1& -\\
 -& 1& -& 1& -& -& -& -& 1& -& -& -& 1& -& -& -& 1& -& 1& -& 1& 1& -& 1& -& -& -& -& -& 1& 1& 1& 1& -& 1& 1\\
 1& -& 1& 1& -& 1& -& 1& 1& 1& -& 1& -& -& 1& -& -& -& -& -& 1& 1& -& -& -& 1& -& 1& -& 1& 1& -& -& -& -& -\\
 -& 1& -& -& 1& -& -& 1& -& 1& -& -& -& -& -& -& -& 1& -& -& -& -& 1& 1& -& 1& 1& 1& 1& 1& 1& 1& -& -& 1& -\\
 -& 1& -& -& 1& -& -& -& 1& 1& 1& 1& -& 1& 1& -& 1& -& -& 1& 1& -& 1& 1& -& -& -& -& -& 1& -& -& -& 1& -& -\\
 -& 1& -& -& -& 1& 1& -& -& -& 1& -& -& -& 1& -& -& -& 1& 1& -& -& -& -& 1& 1& -& 1& -& 1& 1& -& -& 1& 1& 1\\
 -& 1& -& -& -& 1& -& 1& -& 1& -& -& 1& 1& 1& 1& 1& -& -& -& -& 1& 1& -& 1& 1& -& -& -& -& -& -& 1& -& 1& -\\
 1& -& 1& 1& 1& -& 1& 1& -& -& -& -& -& 1& 1& -& 1& -& 1& -& -& -& 1& -& -& -& 1& -& -& 1& -& -& -& -& 1& 1\\
 1& -& -& -& -& 1& -& -& 1& 1& 1& 1& 1& -& -& 1& -& 1& 1& -& -& -& 1& -& -& -& -& -& 1& 1& -& -& -& -& 1& 1\\
 1& -& -& -& -& 1& 1& -& -& -& 1& -& -& -& 1& -& -& -& -& -& 1& 1& 1& 1& 1& -& 1& -& 1& 1& -& 1& 1& -& -& -\\
 -& 1& 1& 1& 1& -& 1& -& 1& -& 1& 1& -& -& -& -& -& 1& -& -& -& 1& 1& -& 1& -& -& 1& -& -& -& -& 1& -& 1& -\\
 -& 1& 1& 1& -& 1& 1& 1& -& -& -& -& 1& -& -& 1& -& 1& -& 1& 1& -& 1& -& -& -& -& -& -& 1& -& 1& -& 1& -& -\\
 1& -& -& -& 1& -& 1& -& -& -& 1& -& 1& 1& -& 1& 1& 1& -& -& 1& -& -& -& -& 1& -& 1& -& 1& 1& -& 1& -& -& -\\
 1& -& -& -& 1& -& -& 1& -& 1& -& -& -& -& -& -& -& 1& 1& 1& 1& 1& 1& -& 1& -& -& -& -& -& 1& -& 1& 1& -& 1\\
 1& -& -& 1& -& -& -& -& 1& -& -& -& 1& -& -& -& 1& -& -& 1& 1& -& 1& -& 1& 1& 1& 1& 1& -& -& -& -& 1& 1& -\\
 -& 1& 1& -& 1& 1& -& 1& 1& -& 1& -& 1& 1& -& -& -& -& -& -& 1& -& -& -& 1& -& 1& -& 1& -& 1& -& -& -& -& 1\\
 1& -& -& 1& -& -& -& 1& -& -& 1& 1& 1& 1& 1& -& -& 1& -& -& -& -& -& 1& 1& -& -& -& -& -& 1& 1& -& 1& 1& -\\ \end{tabular} \right)\] 
\end{scriptsize}

Furthermore, such Hadamard 2-design can be obtained from the tactical decomposition matrix $A_1=[\rho_{ij}]$ for the case with 5 fixed points and 10 orbits of length 3. It has been shown that the matrix $A_1$ expands to the design with the largest automorphism group among constructed structures. When automorphism of order $3$ fixes $2$ points and blocks there is also a TDM leading to the designs with the same automorphism group order, but for the case $F=8$ there is not. This third case is not so rich in automorphism group orders, having $36$ as the largest order appearing. The matrix $A_1$, with the principal submatrix representing orbits being symmetric, is listed as follows.

\begin{scriptsize}
\[ A_1= \left( \begin{tabular}
{c@{\hspace{0.5em}}c@{\hspace{0.5em}}c@{\hspace{0.5em}}c@{\hspace{0.5em}}c@{\hspace{0.5em}}c@{\hspace{0.5em}}c@{\hspace{0.5em}}c@{\hspace{0.5em}}c@{\hspace{0.5em}}c@{\hspace{0.5em}}c@{\hspace{0.5em}}c@{\hspace{0.5em}}c@{\hspace{0.5em}}c@{\hspace{0.5em}}c@{\hspace{0.5em}}}
1&	1&	1&	1&	1&	3&	3&	3&	3&	0&	0&	0&	0&	0&	0\\
1&	1&	1&	1&	1&	3&	0&	0&	0&	3&	3&	3&	0&	0&	0\\
1&	1&	1&	1&	1&	0&	3&	0&	0&	3&	0&	0&	3&	3&	0\\
1&	1&	1&	1&	1&	0&	0&	3&	0&	0&	3&	0&	3&	0&	3\\
1&	1&	1&	1&	1&	0&	0&	0&	3&	0&	0&	3&	0&	3&	3\\
1&	1&	0&	0&	0&	3&	1&	1&	1&	1&	1&	1&	2&	2&	2\\
1&	0&	1&	0&	0&	1&	3&	1&	1&	1&	2&	2&	1&	1&	2\\
1&	0&	0&	1&	0&	1&	1&	3&	1&	2&	1&	2&	1&	2&	1\\
1&	0&	0&	0&	1&	1&	1&	1&	3&	2&	2&	1&	2&	1&	1\\
0&	1&	1&	0&	0&	1&	1&	2&	2&	3&	1&	1&	1&	1&	2\\
0&	1&	0&	1&	0&	1&	2&	1&	2&	1&	3&	1&	1&	2&	1\\
0&	1&	0&	0&	1&	1&	2&	2&	1&	1&	1&	3&	2&	1&	1\\
0&	0&	1&	1&	0&	2&	1&	1&	2&	1&	1&	2&	3&	1&	1\\
0&	0&	1&	0&	1&	2&	1&	2&	1&	1&	2&	1&	1&	3&	1\\
0&	0&	0&	1&	1&	2&	2&	1&	1&	2&	1&	1&	1&	1&	3\\
\end{tabular} \right)\] 
\end{scriptsize}

\subsection{Results for Hadamard 3-designs with parameters $3$-$(36,18,8)$}
Although it is possible to construct 3-designs directly from adequate tactical decomposition matrices we take an advantage of already constructed Hadamard matrices. According to previously described procedure, every Hadamard 2-design lead to the one or more Hadamard 3-designs. It has shown that there are 97662, 1726 and 3413 $3$-$(36,18,8)$ designs arising from related Hadamard matrices of Cyc* type (those obdained from Hadamard 2-designs admitting an automorphism of order 3), with 2, 5 and 8 fixed points, respectively. There are 102722 non-isomorphic designs among these three classes. Contrary to the results for $2$-$(35,17,8)$, here we get structures without any other automorphism but identity.  It follows the full list of appearing automorphism group orders and related frequencies: 1 (77909), 2 (931), 3 (22734), 4 (130), 6 (778), 8 (22), 12 (103), 16 (6), 18 (31), 21 (1), 24 (28), 27 (2), 32 (2), 36 (7), 48 (3), 54 (4), 60 (1), 64 (2), 72 (4), 96 (5), 108 (1), 144 (2),  168 (1), 192 (1), 384 (1), 420 (1), 432 (1), 576 (1), 720 (1), 40320 (1).

For the ACyc case of the construction there are 195543, 13317 and 65169 non-isomorphic designs with parameters $3$-$(36,18,8)$ arising from related 2-designs. Finally, there are 272116 non-isomorphic designs with these parameters among these three sets of designs. The summarized result is presented in the next proposition.

\begin{prop}
There are at least $272116$ Hadamard $3$-designs with parameters $(36,18,8)$, where $237444$ of them do not admit any non-trivial automorphism.
\end{prop}

\begin{center}
\begin{figure}[h!]
\begin{scriptsize}
\[ \begin{tabular}{|c@{\hspace{0.4em}}c@{\hspace{0.4em}}c@{\hspace{0.4em}}c@{\hspace{0.4em}}c@{\hspace{0.4em}}c@{\hspace{0.4em}}c@{\hspace{0.4em}}c@{\hspace{0.4em}}c@{\hspace{0.4em}}c@{\hspace{0.4em}}|c@{\hspace{0.4em}}c@{\hspace{0.4em}}c@{\hspace{0.4em}}|c@{\hspace{0.4em}}c@{\hspace{0.4em}}c@{\hspace{0.4em}}|c@{\hspace{0.4em}}c@{\hspace{0.4em}}c@{\hspace{0.4em}}|c@{\hspace{0.4em}}c@{\hspace{0.4em}}c@{\hspace{0.4em}}|c@{\hspace{0.4em}}c@{\hspace{0.4em}}c@{\hspace{0.4em}}|c@{\hspace{0.4em}}c@{\hspace{0.4em}}c@{\hspace{0.4em}}|c@{\hspace{0.4em}}c@{\hspace{0.4em}}c@{\hspace{0.4em}}|c@{\hspace{0.4em}}c@{\hspace{0.4em}}c@{\hspace{0.4em}}|c@{\hspace{0.4em}}c@{\hspace{0.4em}}c@{\hspace{0.4em}}|c@{\hspace{0.4em}}c@{\hspace{0.4em}}c@{\hspace{0.4em}}|} 
 \hline
1&	1&	1&	1&	1&	&	&	&	&	&	1&	1&	1&	1&	1&	1&	1&	1&	1&	1&	1&	1&	1&	1&	1&	1&	1&	1&	1&	1&	1&	1&	1&	1&	1&	1&	1&	1&	1&	1\\
1&	1&	1&	1&	1&	&	&	&	&	&	1&	1&	1&	1&	1&	1&	1&	1&	1&	1&	1&	1&	1&	1&	1&	1&	1&	1&	1&	1&	1&	1&	1&	1&	1&	1&	1&	1&	1&	1\\
1&	1&	1&	1&	1&	&	&	&	&	&	1&	1&	1&	1&	1&	1&	1&	1&	1&	1&	1&	1&	1&	1&	1&	1&	1&	1&	1&	1&	1&	1&	1&	1&	1&	1&	1&	1&	1&	1\\
1&	1&	1&	1&	1&	&	&	&	&	&	&	&	&	&	&	&	&	&	&	&	&	&	&	&	&	&	&	&	&	&	&	&	&	&	&	&	&	&	&	\\
1&	1&	1&	1&	1&	&	&	&	&	&	&	&	&	&	&	&	&	&	&	&	&	&	&	&	&	&	&	&	&	&	&	&	&	&	&	&	&	&	&	\\
1&	1&	1&	1&	1&	&	&	&	&	&	&	&	&	&	&	&	&	&	&	&	&	&	&	&	&	&	&	&	&	&	&	&	&	&	&	&	&	&	&	\\ \hline
1&	1&	&	&	&	1&	1&	1&	&	&	1&	1&	1&	1&	&	&	1&	&	&	1&	&	&	&	1&	1&	&	1&	1&	&	1&	1&	1&	&	1&	1&	&	1&	1&	&	1\\
1&	1&	&	&	&	1&	1&	1&	&	&	1&	1&	1&	&	1&	&	&	1&	&	&	1&	&	1&	&	1&	1&	&	1&	1&	&	1&	1&	1&	&	1&	1&	&	1&	1&	\\
1&	1&	&	&	&	1&	1&	1&	&	&	1&	1&	1&	&	&	1&	&	&	1&	&	&	1&	1&	1&	&	1&	1&	&	1&	1&	&	&	1&	1&	&	1&	1&	&	1&	1\\ \hline
1&	&	1&	&	&	&	1&	1&	&	1&	1&	&	&	1&	1&	1&	1&	&	&	1&	&	&	&	1&	1&	1&	&	&	1&	&	&	&	1&	&	&	1&	&	1&	&	1\\
1&	&	1&	&	&	&	1&	1&	&	1&	&	1&	&	1&	1&	1&	&	1&	&	&	1&	&	1&	&	1&	&	1&	&	&	1&	&	&	&	1&	&	&	1&	1&	1&	\\
1&	&	1&	&	&	&	1&	1&	&	1&	&	&	1&	1&	1&	1&	&	&	1&	&	&	1&	1&	1&	&	&	&	1&	&	&	1&	1&	&	&	1&	&	&	&	1&	1\\ \hline
1&	&	&	1&	&	1&	&	1&	&	1&	1&	&	&	1&	&	&	1&	1&	1&	1&	&	&	1&	&	&	&	1&	1&	1&	&	&	&	1&	&	1&	&	1&	&	1&	\\
1&	&	&	1&	&	1&	&	1&	&	1&	&	1&	&	&	1&	&	1&	1&	1&	&	1&	&	&	1&	&	1&	&	1&	&	1&	&	&	&	1&	1&	1&	&	&	&	1\\
1&	&	&	1&	&	1&	&	1&	&	1&	&	&	1&	&	&	1&	1&	1&	1&	&	&	1&	&	&	1&	1&	1&	&	&	&	1&	1&	&	&	&	1&	1&	1&	&	\\ \hline
1&	&	&	&	1&	1&	1&	&	&	1&	1&	&	&	1&	&	&	1&	&	&	1&	1&	1&	1&	&	&	1&	&	&	&	1&	1&	1&	&	1&	&	1&	&	&	1&	\\
1&	&	&	&	1&	1&	1&	&	&	1&	&	1&	&	&	1&	&	&	1&	&	1&	1&	1&	&	1&	&	&	1&	&	1&	&	1&	1&	1&	&	&	&	1&	&	&	1\\
1&	&	&	&	1&	1&	1&	&	&	1&	&	&	1&	&	&	1&	&	&	1&	1&	1&	1&	&	&	1&	&	&	1&	1&	1&	&	&	1&	1&	1&	&	&	1&	&	\\ \hline
&	1&	1&	&	&	&	1&	1&	1&	&	1&	&	&	1&	&	&	&	1&	1&	&	1&	1&	&	&	&	&	1&	1&	&	1&	1&	&	1&	&	&	1&	&	1&	&	1\\
&	1&	1&	&	&	&	1&	1&	1&	&	&	1&	&	&	1&	&	1&	&	1&	1&	&	1&	&	&	&	1&	&	1&	1&	&	1&	&	&	1&	&	&	1&	1&	1&	\\
&	1&	1&	&	&	&	1&	1&	1&	&	&	&	1&	&	&	1&	1&	1&	&	1&	1&	&	&	&	&	1&	1&	&	1&	1&	&	1&	&	&	1&	&	&	&	1&	1\\ \hline
&	1&	&	1&	&	1&	&	1&	1&	&	1&	&	&	&	1&	1&	1&	&	&	&	1&	1&	&	1&	1&	&	&	&	&	1&	1&	&	1&	&	1&	&	1&	&	1&	\\
&	1&	&	1&	&	1&	&	1&	1&	&	&	1&	&	1&	&	1&	&	1&	&	1&	&	1&	1&	&	1&	&	&	&	1&	&	1&	&	&	1&	1&	1&	&	&	&	1\\
&	1&	&	1&	&	1&	&	1&	1&	&	&	&	1&	1&	1&	&	&	&	1&	1&	1&	&	1&	1&	&	&	&	&	1&	1&	&	1&	&	&	&	1&	1&	1&	&	\\ \hline
&	1&	&	&	1&	1&	1&	&	1&	&	1&	&	&	&	1&	1&	&	1&	1&	1&	&	&	&	1&	1&	&	1&	1&	&	&	&	1&	&	1&	&	1&	&	&	1&	\\
&	1&	&	&	1&	1&	1&	&	1&	&	&	1&	&	1&	&	1&	1&	&	1&	&	1&	&	1&	&	1&	1&	&	1&	&	&	&	1&	1&	&	&	&	1&	&	&	1\\
&	1&	&	&	1&	1&	1&	&	1&	&	&	&	1&	1&	1&	&	1&	1&	&	&	&	1&	1&	1&	&	1&	1&	&	&	&	&	&	1&	1&	1&	&	&	1&	&	\\ \hline
&	&	1&	1&	&	&	&	1&	1&	1&	1&	&	1&	&	1&	&	&	1&	&	1&	&	1&	1&	&	1&	1&	&	1&	&	1&	&	1&	1&	1&	&	&	1&	&	&	1\\
&	&	1&	1&	&	&	&	1&	1&	1&	1&	1&	&	&	&	1&	&	&	1&	1&	1&	&	1&	1&	&	1&	1&	&	&	&	1&	1&	1&	1&	1&	&	&	1&	&	\\
&	&	1&	1&	&	&	&	1&	1&	1&	&	1&	1&	1&	&	&	1&	&	&	&	1&	1&	&	1&	1&	&	1&	1&	1&	&	&	1&	1&	1&	&	1&	&	&	1&	\\ \hline
&	&	1&	&	1&	&	1&	&	1&	1&	1&	&	1&	&	1&	&	1&	&	1&	&	1&	&	1&	&	1&	&	1&	&	1&	&	1&	&	&	1&	1&	1&	1&	&	&	1\\
&	&	1&	&	1&	&	1&	&	1&	1&	1&	1&	&	&	&	1&	1&	1&	&	&	&	1&	1&	1&	&	&	&	1&	1&	1&	&	1&	&	&	1&	1&	1&	1&	&	\\
&	&	1&	&	1&	&	1&	&	1&	1&	&	1&	1&	1&	&	&	&	1&	1&	1&	&	&	&	1&	1&	1&	&	&	&	1&	1&	&	1&	&	1&	1&	1&	&	1&	\\ \hline
&	&	&	1&	1&	1&	&	&	1&	1&	1&	&	1&	1&	&	1&	&	1&	&	&	1&	&	&	1&	&	1&	&	1&	1&	&	1&	&	&	1&	&	&	1&	1&	1&	1\\
&	&	&	1&	1&	1&	&	&	1&	1&	1&	1&	&	1&	1&	&	&	&	1&	&	&	1&	&	&	1&	1&	1&	&	1&	1&	&	1&	&	&	1&	&	&	1&	1&	1\\
&	&	&	1&	1&	1&	&	&	1&	1&	&	1&	1&	&	1&	1&	1&	&	&	1&	&	&	1&	&	&	&	1&	1&	&	1&	1&	&	1&	&	&	1&	&	1&	1&	1\\
 \hline \end{tabular} \] 
\end{scriptsize}
\caption{Submatrix of an incidence matrix of Hadamard 3-design ${\cal D}_1$ with 30 points and 70 lines, having $|Aut({\cal D}_1)|=40320$ (with 0's omitted). It is seen 40 out of 70 lines, while the rest of lines are the complement of columns 11-40.} \label{Fig1}
\end{figure}
\end{center}

All these constructed Hadamard 3-designs have the next automorphism group orders: 1 (237444), 2 (9304), 3 (22734), 4 (1462), 6 (778), 8 (127), 9 (8), 10 (2), 12 (103), 16 (46), 18 (31), 21 (1), 24 (28), 27 (2), 32 (6), 36 (7), 48 (3), 54 (4), 60 (1), 64 (5), 72 (4), 96 (5), 108 (1), 144 (2), 168 (1), 192 (1), 384 (1), 420 (1), 432 (1), 576 (1), 720 (1), 40320 (1).

The vast majority of designs do not admit any non-trivial automorphism, while there is one design with the largest automorphism group, $|Aut({\cal D})|=40320$. The only automorphism group order appearing within this construction and not present for Cyc construction is 10. All other automorphism group orders are the same with possibly different frequencies. 

As an example let mention that the matrix $H_{(40320)}$, that arises from the most symmetric design amoung constructed Hadamard 2-designs, gives 36 incidence matrices of the design with parameters $3$-$(36,18,8)$. Up to isomorphism, these 36 matrices results with one design ${\cal D}_1$, again with the automorphism group order of 40320. Denote one of these incidence matrices by $M$. Having in mind that an automorphism of order 3 form a tactical decomoposition (while in general do not every TDM lead to an incidence matrix of a design), there must be a matrix equivalent to $M$ that consist only of cyclic matrices of order 3 in addition to all-1 and all-0 matrices. More precisely, ${\cal D}_1$ is represented by an incidence matrix having 200 cyclic matrices of order 3. It has shown that half of these matrices is the complement of the rest of matrices, as it is depicted in Figure \ref{Fig1} (0's are omitted).

It is worth mentioning that it is possible to construct self-dual codes from Hadamard matrices (for details see \cite{Slo}). In particular, there are problems of interest in coding theory concerning the existence of an extremal self-dual $[72,36,16]$ code \cite{SpTo}. Presented results are tidy related to this issue since a code with this parameters can be constructed from Hadamard matrices of order 36 with automorphism of order 1,2,3,5 or 7 \cite{DZD}.

Table \ref{Table1} presents previously described results for designs with parameters $2$-$(35,17,8)$, $3$-$(36,18,8)$ and Hadamard matrices of order 36.

\begin{center}
\begin{table}[h]
$$
\begin{array}{rrrrrrr}
\hline
&\multicolumn{2}{c}{2$-$(35,17,8)}& \multicolumn{2}{c} {HM(36)} &  \multicolumn{2}{c}{3$-$(36,18,8)}\\
\hline
F & \textrm{ Cyc}&\textrm{ ACyc}& \textrm{ Cyc*}& \textrm{ ACyc*}&\textrm{ Cyc*}&\textrm{ ACyc*}\\ \hline
2 &63635 & 89686 & 7238 &10174  & 97662&195543\\
5 & 3698& 17627&158&635 &1726&13317\\
8 & 14692&  155550&259&2176&3413&65169\\ \hline
\textrm{all} & 81973& 262511&7631&12833&102722&272116\\ \hline
\end{array}
$$
\caption{The number of constructed non-isomorphic Hadamard 2-designs with 35 points, inequivalent Hadamard matrices of order 36 and the number of non-isomorphic Hadamard 3-designs with 36 points and 70 lines. \label{Table1}}
\end{table}
\end{center}




\begin{thebibliography}{00}

\bibitem{BFW}
I. Bouyukliev, V. Fack, J. Winne, 2-(31,15,7), 2-(35,17,8) and 2-(36,15,6) designs with automorphisms of odd prime order, and their related Hadamard matrices and codes, Des.Codes Cryptogr., 51 (2009), 105-122.

\bibitem{BBI}
L. Berardi, M. Buratti, S. Innamorati, 4-Blocked Hadamard 3-designs, Disc. Math., 174 (1997), 35-46.

\bibitem{CEPU}
V. \'Cepuli\'c, On symmetric block designs $(40,13,4)$ with automorphisms of order $5$, Discrete Math., 128 (1994), 45-60.

\bibitem{DZD}
R. Dontcheva, A.J. Zanten, S. M. Dodunekov, Binary self-dual codes with automorphism of composite order, IEEE Trans. Inform. Theory, 50/2 (2004), 311-318.


\bibitem{HuPu}
D.R. Hughes, F.C. Piper, Design theory, Cambridge University Press, Cambridge, 1985.

\bibitem{ITT}
Y. J. Ionin and Tran van Trung, Symmetric Designs, in CRC Handbook of Combinatorial Designs Second Edition, C. J. Colbourn and J. H. Dinitz (Editors), CRC Press, Boca Raton, FL, (2007) 110-124.

\bibitem{Jank}
Z. Janko, The Existence of a Bush-type Hadamard Matrix of Order 36 and Two New Infinite Classes of Symmetric Designs, J. Comb. Theory, Series A, 95 (2001), 360-364.

\bibitem{JTT}
Z. Janko and Tran van Trung,
Construction of a new symmetric block design for $(78,22,6)$ with the help of tactical decompositions, J. Comb. Theory Ser. A, { 40} (1985), 451-455.

\bibitem{KaOs} 
P. Kaski and P. R. J. \"Osterg\r{a}rd, Classification Algorithms for Codes and Designs, Springer, Berlin, 2006.

\bibitem{Khar}
H. Kharaghani, B. Tayfeh-Rezaie, On the classification of Hadamard matrices of order 32, J. Combin. Des., 18 (2010), 328–336.

\bibitem{Krc}
V. Kr\v cadinac, Steiner $2$-designs $S(2,5,41)$ with automorphisms of order $3$, Journal of Combinatorial Mathematics and Combinatorial Computing, {43} (2002), 83-99.

\bibitem{MaPa}
I. Martinjak, M.O.Pav\v cevi\'c, Symmetric designs possesing tactical decompositions, Advances in Mathematics of Communications, 10 (2011), 199-208.

\bibitem{McK1}
B. D. McKay, Hadamard equivalence via graph isomorphism, Disc. Math., 27 (1979), 213-214.

\bibitem{McK2} 
B. D. McKay, {\sf nauty} user's guide (version 1.5), Technical Report TR-CS-90-02, Department of Computer Science, Australian National University, 1990.

\bibitem{Slo}
N.J.A. Sloane, Is there a (72,36), de=16 self-dual code?, IEEE Trans. Inform. Theory, 19 (1973), 251.

\bibitem{SpTo}
E. Spence, V.D. Tonchev, Extremal self-dual codes from symmetric designs, Disc. Math., 110 (1992), 265-268.

\bibitem{Tonc}
V. D. Tonchev, Hadamard Matrices of Order 36 with Automorphisms of Order 17, Nagoya Math. J., 104 (1986), 163-174.



\end{thebibliography}




\end{document}